\documentclass[11pt]{article}
\usepackage{mathtools}
\usepackage{amsfonts}
\usepackage{amsmath}
\usepackage{amssymb}
\usepackage{amsthm}
\usepackage{xfrac}
\usepackage{graphicx}
\usepackage{todonotes}
\usepackage{comment}
\usepackage{subfigure}
\usepackage{theoremref}
\usepackage{authblk}
\usepackage{citesort}
\usepackage{multirow}

\usepackage{diagbox}

\newtheorem{theorem}{Theorem}

\newtheorem{remark}{Remark}

\newcommand{\norm}[1]{\left\lVert#1\right\rVert}

\def \ol#1{\overline{#1}}      

\begin{document}
\title{Preconditioned Iterative Methods for Diffusion Problems\\ with High-Contrast Inclusions}
\author[1]{Yuliya Gorb\thanks{gorb@math.uh.edu, corresponding author}}
\author[2]{{Vasiliy Kramarenko\thanks{kramarenko.vasiliy@gmail.com}}}
\author[1]{Yuri Kuznetsov\thanks{kuz@math.uh.edu}}
\affil[1]{Department of Mathematics,  University of Houston, Houston, TX 77204}
\affil[2]{
Marchuk Institute of Numerical Mathematics of the Russian Academy of Sciences,
Moscow, Russian Federation}
\date{}

\maketitle

\begin{abstract}
\noindent 
This paper concerns robust numerical treatment of an elliptic PDE with high contrast coefficients, for which classical finite-element discretizations yield ill-conditioned linear systems. This paper introduces a procedure by which the discrete system obtained from a linear finite element discretization of the given continuum problem is converted into an equivalent linear system of the saddle point type.  Then three preconditioned iterative procedures -- preconditioned Uzawa, preconditioned Lanczos, and PCG for the square of the matrix -- are discussed for a special type of the application, namely, highly conducting particles distributed in the domain. 
Robust preconditioners for solving the derived saddle point problem are proposed and investigated. Robustness with respect to the contrast parameter and the discretization scale is also justified. Numerical examples support theoretical results and demonstrate independence of the number of iterations of the proposed iterative schemes on the contrast in parameters of the problem and the mesh size. 
\end{abstract}


\section{Introduction}

\begin{figure}[h!tb]
	\centering
	\includegraphics[width=0.5\textwidth]{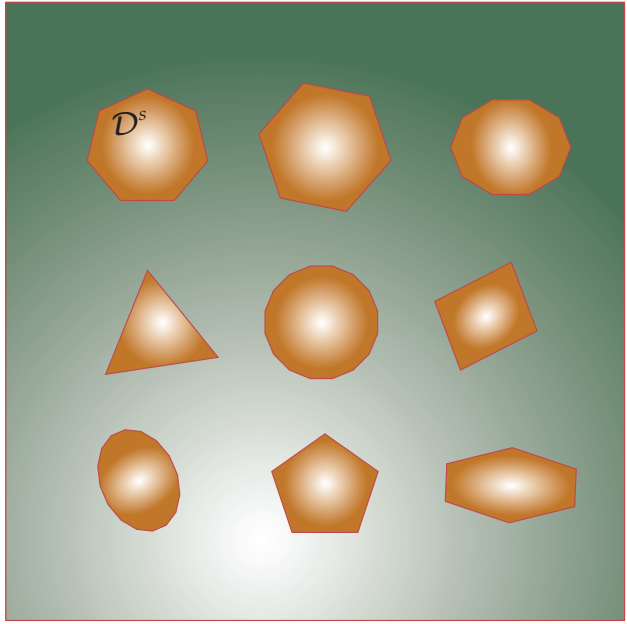}
	\caption{An example of $\mathcal{D}^s$.}
	\label{f:kuz_fig_ex0}
\end{figure}

In this paper, we consider iterative solutions of the linear system arising from the discretization of a diffusion problem
\begin{equation} \label{E:problem-intro}
-\nabla \cdot \left[  \sigma(x) \nabla u \right]  = f, \quad x\in  \Omega
\end{equation}
with appropriate boundary conditions on $\Gamma = \partial \Omega$. Below, in our theoretical consideration and numerical tests, we will assume the homogeneous Dirichlet boundary conditions on $\Gamma$.
The main focus of this work is on the case when the coefficient function $ \sigma(x) \in L^{\infty} (\Omega)$ varies largely within the domain $\Omega$, that is,
\[
\kappa = \frac{\sup_{x\in \Omega}  \sigma(x)}{\inf_{x\in \Omega}  \sigma(x)} \gg 1.
\]
We assume that $\Omega$ is a bounded domain $\Omega \subset \mathbb{R}^2$, that contains $m\geq 1$ disjoint polygonal 
subdomains $\mathcal{D}^s$, $s \in \{1,\ldots,m\}$, see Figure \ref{f:kuz_fig_ex0}, in which $\sigma$ is ``large'', e.g. of order $O(\kappa)$, but remains of $O(1)$ in the domain outside of $\mathcal{D} := \cup_{s=1}^m \mathcal{D}^s$. 

A P1-FEM discretization of this problem results in a linear system
\begin{equation} \label{E:linsys}
\bold{A}_{\sigma} \,\overline{u}  = \overline{f},
\end{equation}
with a large, sparse, symmetric and positive definite (SPD) matrix $\bold{A}_{\sigma}$. 
A major issue in numerical treatments of \eqref{E:problem-intro} with the coefficient $\sigma$ discussed above, is that the high contrast leads to an ill-conditioned matrix $\bold{A}_{\sigma}$. Indeed, if $h$ is the discretization scale, then the condition number of the resulting stiffness matrix $\bold{A}_{\sigma}$ grows proportionally to $h^{-2}$ with the coefficient of proportionality linearly depending on $\kappa$. Because of that, the high contrast problems have been a subject of an active research recently, see e.g. \cite{ah02,agks08,gk18,bgw14}.

Our main goal here is robust numerical treatment of the described problem. For that, we introduce an additional variable that allows us to replace \eqref{E:linsys} with an equivalent formulation in terms of a linear system 
 \begin{equation} \label{E:LS}
\boldsymbol{\mathcal{A}} \, \overline{x}=\overline{F}, \quad \mbox{with} \quad 
\overline{F}=  \begin{bmatrix} \overline{ \mathrm{f}} \\ \overline{0} \end{bmatrix} ,
 \end{equation}
 and a {\it saddle point matrix}  $\boldsymbol{\mathcal{A}}$ written in the block form:
 \begin{equation} \label{E:Ablock}
\boldsymbol{\mathcal{A}}=
\begin{bmatrix}
\bold{A} & \bold{B}^T \\
\bold{B} & -\bold{C}
\end{bmatrix},
 \end{equation}
where $\bold{A} \in \mathbb{R}^{N \times N}$ is SPD, $\bold{B}  \in \mathbb{R}^{n \times N} $ is
rank deficient, and $\bold{C} \in \mathbb{R}^{n \times n} $ is an SPD matrix. 
Below, we discuss three iterative procedures -- preconditioned Uzawa (PU) method for the system with an SPD Schur complement matrix; preconditioned Lanczos (PL) method for solving \eqref{E:LS}; and preconditioned conjugate gradient (PCG) method for an equivalent system with an SPD matrix. Then we propose a robust block-diagonal preconditioner 
\[
\boldsymbol{\mathcal{H}}=
\begin{bmatrix}
\mathcal{H}_{\mathrm{A}} & 0\\
0 & \mathcal{H}_{\mathrm{S}}
\end{bmatrix}, 
\] 
for solving \eqref{E:LS}-\eqref{E:Ablock} with these three  iterative methods. The main feature of the proposed preconditioners is that convergence rates of discussed iterative schemes are independent of the contrast parameter $\kappa \gg 1$ and the discretization size $h>0$.  A rigorous justification of the latter statement is based on the evaluation of the eigenvalues of the matrix
$\boldsymbol{\mathcal{H}}\boldsymbol{\mathcal{A}}$, which are proven to be in the union of two intervals  $[\mu_-^1,\mu_-^2]\cup[\mu_+^1,\mu_+^2]$, where $\mu_-^1< \mu_-^2 < 0 < \mu_+^1 < \mu_+^2$. Assuming that the mesh on $\Omega$ is regularly-shaped and quasi-uniform, we demonstrate that constants $\mu_\pm^i$ ($i=1,2$) are independent of the discretization scale $h$ and the number of inclusions. If, in addition, we assume that particles are located at distances comparable to their sizes, then $\mu_\pm^i$ ($i=1,2$) are independent of the diameters of $\mathcal{D}^s$, $s \in \{1,\ldots,m\}$, their locations, and distances between them.
The numerical experiments on simple test cases support  theoretical findings and demonstrate independence of convergence  rates of the proposed iterative schemes on parameters indicated above. These numerical tests are performed for a two-dimensional problem, whereas  theoretical results remain true for  three dimensions as well.

The development of efficient preconditioners for saddle point problems has been an active area of research since early 1990s, see e.g. \cite{bpx90,bpx97,rus90,wath,erv90}. The main feature of the problem considered in this paper is that we deal with a special type of saddle point matrices that, in particular, contains a rank deficient block $\bold{B}$.
Also, this paper proposes a very special form for the block $ \mathcal{H}_{\mathrm{S}}$ of $\boldsymbol{\mathcal{H}}$, see \eqref{E:Uzawa-prec} in Section \ref{S:methods}, utilized in three methods that yields theoretical results mentioned above.
Moreover,  one of the iterative procedures that we employed in this paper is {\it Lanczos method} \cite{lanc50,paige} that, as would be evident from our numerical experiments below, has demonstrated significant advantages over the other methods with respect to the arithmetic cost.

Finally, we point out that robust numerical treatment of the described problem
is crucial in developing the mutiscale strategies for models of composite materials with 
highly conducting particles. The latter find their application in particulate flows, subsurface flows in natural porous formations,
electrical conduction in composite materials, and medical and geophysical imaging.

\vspace{5pt}

The paper is organized as follows. In Section \ref{S:form}, the mathematical problem formulation is presented including the derivation of the saddle point problem of the type \eqref{E:LS}-\eqref{E:Ablock}. 
Section \ref{S:methods} discusses three iterative methods (preconditioned Uzawa, preconditioned Lanczos and PCG, mentioned above) for solving system \eqref{E:LS}-\eqref{E:Ablock} and proposes efficient preconditioners for all of them. The main theoretical results, which are the estimates for the eigenvalues of the matrix $\boldsymbol{\mathcal{H}}\boldsymbol{\mathcal{A}}$, are stated and proven in Section \ref{S:eig_est}. Numerical experiments based on simple test scenarios are presented in Section \ref{S:numerics}. Conclusions are discussed in Section \ref{S:concl}.

\vspace{8pt}

\noindent {\bf Acknowledgements.} 
Y. Gorb has been supported by the NSF grant DMS-$1350248$.

\section{Problem Formulation} \label{S:form}

\subsection{Equivalent variational formulations}

Consider an open, a bounded domain $\Omega \subset \mathbb{R}^2$ with a piece-wise smooth boundary $\Gamma:=\partial \Omega$, that contains $m\geq 1$ subdomains $\mathcal{D}^s$ with piece-wise smooth boundaries $\Gamma_s:=\partial \Omega_s$, $s\in\{1,\ldots,m\}$, see Figure \ref{f:kuz_fig_ex0}. 
Assume that  $\Gamma_s \cap \Gamma_t = \emptyset$ when $s\neq t$, and $\Gamma \cap \Gamma_s= \emptyset $, $s\in\{1,\ldots,m\}$. 
For simplicity, we assume that $\Omega$ and $\mathcal{D}^s$ are polygons. 
The union of $\mathcal{D}^s$ is denoted by $\mathcal{D}$.
In the domain $\Omega$, we consider the following elliptic problem
\begin{equation} \label{E:pde-form}
\left\{
\begin{array}{r l l}
-\nabla \cdot \left[  \sigma(x) \nabla u \right] & = f, & x\in  \Omega  \\[2pt]
u & = 0, & x\in  \Gamma
\end{array}
\right.
\end{equation}
with the source term $f \in L^2(\Omega)$, and the coefficient $\sigma$ that varies largely inside the domain $\Omega$. 
In this paper, we are focused on the case when $\sigma$ is a piecewise constant function given by
\begin{equation} \label{E:sigma}
\sigma(x)=
\begin{cases}
1, & x\in  \Omega\setminus \overline{\mathcal{D}}\\
\displaystyle 1+\frac{1}{\varepsilon_s}, &  x\in \mathcal{D}^s, ~s\in\{1,\ldots,m\}
\end{cases}
\end{equation}
with $\displaystyle 0< \varepsilon_s\equiv \mbox{const} \leqslant 1$, $s\in\{1,\ldots,m\}$. The standard variational formulation of \eqref{E:pde-form} is 
\begin{equation} \label{E:var-form-orig}
\mbox{Find}~ u \in V:=H^1_0(\Omega) ~ 
\mbox{such that}~
\int_{\Omega} \nabla u \cdot  \nabla v ~dx +\sum_{s=1}^{m}\frac{1}{\varepsilon_s} \int_{\mathcal{D}^s} \nabla u \cdot  \nabla v ~dx = \int_{\Omega} f v ~dx, ~ \forall v\in V.
\end{equation}
We introduce new variables $p_s \in H^{1}(\mathcal{D}^s)$ via
\begin{equation} \label{E:p}
p_s =\frac{1}{\varepsilon_s} u_s +c_s  \quad \mbox{in} \quad\mathcal{D}^s, \quad  s\in \{1,\ldots,m\}, \quad \mbox{where} ~u_s = u|_{ \mathcal{D}^s},
\end{equation}
and $c_s$ are arbitrary constants, $s\in\{1,\ldots,m\}$.
With that, we replace formulation \eqref{E:var-form-orig} with the new one, namely,
\[
\mbox{Find}~ u \in V~  \mbox{and} ~ p_s \in  \left.V_s:=H^1(\mathcal{D}^s)=V\right|_{\mathcal{D}^s}, ~s\in \{1,\ldots,m\}, ~ \mbox{such that}  
\]
\begin{equation} \label{E:var-form-new1}
\int_{\Omega} \nabla u \cdot  \nabla v ~dx + \sum_{t=1}^m  \int_{\mathcal{D}^t} \nabla p_t \cdot  \nabla v ~dx = \int_{\Omega} f v ~dx, ~ \forall v\in V, \end{equation}
\begin{equation} \label{E:var-form-new2}
 \int_{\mathcal{D}^s} \nabla u \cdot  \nabla w ~dx -   \varepsilon_s \int_{\mathcal{D}^s} \nabla p_s
 \cdot  \nabla w ~dx=0 ,  ~ \forall w \in V_s,~ s\in\{1,\ldots,m\}.
\end{equation}
Two formulations \eqref{E:var-form-orig} and \eqref{E:var-form-new1}-\eqref{E:var-form-new2} are equivalent in the sense that their solutions $u\in H^1(\Omega)$ coincide, and any solution $p_s\in V_s$ of \eqref{E:var-form-new1}-\eqref{E:var-form-new2} is 
equal to the function $\frac{1}{\varepsilon_s} u_s +c_s$  with an appropriate constant $c_s$, $s\in\{1,\ldots,m\}$.
For the uniqueness of $p_s$, we can either demand 
\begin{equation} \label{E:uniq-p}
\int_{\mathcal{D}^s}   p_s ~dx = 0, \quad s\in\{1,\ldots,m\} , 
\end{equation}
or modify the formulation \eqref{E:var-form-new2} as follows
\begin{equation} \label{E:var-form-ours}
\begin{array}{l l}
\displaystyle 
\mbox{Find}~ u \in V~  \mbox{and} ~ p_s \in V_s~ \mbox{such that}   
\int_{\Omega} \nabla u \cdot  \nabla v ~dx + \sum_{t=1}^m  \int_{\mathcal{D}^t} \nabla p_t \cdot  \nabla v ~dx = \int_{\Omega} f v ~dx, ~ \forall v\in V,  \\[7pt]
\displaystyle  \int_{\mathcal{D}^s} \nabla u \cdot  \nabla w ~dx - \int_{\mathcal{D}^s} \nabla p_s \cdot  \nabla w ~dx - \frac{1}{|\mathcal{D}^s|} \left[\int_{\mathcal{D}^s}  p_s ~dx \right]  \left[\int_{\mathcal{D}^s} w ~dx \right]  =0 ,  ~ \forall w \in V_s,~ s\in\{1,\ldots,m\},
\end{array}
\end{equation}
where  $|\mathcal{D}^s|$ is the  area of the particle $\mathcal{D}^s$. It is obvious, that solutions $p_s$, $s\in\{1,\ldots,m\}$, of \eqref{E:var-form-ours} satisfy condition  \eqref{E:uniq-p}, and the above constants $c_s$ are defined by 
\[
c_s=-\frac{1}{\varepsilon_s}  \int_{\mathcal{D}^s} u \,dx,  \quad ~ s\in\{1,\ldots,m\}.
\]

\subsection{Discretization of \eqref{E:var-form-ours} and Description of the Saddle Point Problem}

Let $\Omega_h$ be a 
triangular mesh on $\Omega$.  Assume that $\Omega_h$  is conforming with boundaries $\Gamma$ and $\Gamma_s$, $s\in\{1,\ldots,m\}$, that is,  $\Gamma$ and $\Gamma_s$ are the unions of the triangular sides. 
We define $\left.  \mathcal{D}^s_h = \Omega_h\right|_{\mathcal{D}^s}$, $s\in\{1,\ldots,m\}$, and $ \mathcal{D}_h := \cup_{s=1}^m   \mathcal{D}^s_h$.

We now choose a FEM space $V_h \subset H_0^1 (\Omega)$ to be the space of linear finite-element functions defined on  $\Omega_h$, and $V^s_h:= V_{h}|_{\mathcal{D}_h^s}$,  $s\in\{1,\ldots,m\}$. Then, the FEM discretization \cite{bs08} of \eqref{E:var-form-ours} reads as follows: 
\begin{equation} \label{E:FEM-form}
\begin{array}{l l}
\mbox{Find}\quad u_h \in V_h \quad \mbox{and} \quad p_h = (p^1_h,\ldots,p^m_h) \quad \mbox{with}\quad  p^s_h \in V^s_h  \quad \mbox{such that}\\[2pt]
\displaystyle \int_{\Omega} \nabla u_h \cdot  \nabla v_h ~dx + \int_{\mathcal{D}} \nabla p_h \cdot  \nabla v_h ~dx = \int_{\Omega} f \, v_h ~dx, \quad \forall v_h\in V_h,\\[5pt]
\displaystyle \int_{\mathcal{D}^s} \nabla u_h \cdot  \nabla w^s_h ~dx -  \varepsilon_s \int_{\mathcal{D}^s} \nabla p^s_h  \cdot  \nabla w^s_h ~dx - \frac{1}{|\mathcal{D}^s|}  \left[\int_{\mathcal{D}^s} p^s_h ~dx \right]  \left[\int_{\mathcal{D}^s} w^s_h  ~dx \right]  =0 , ~ \forall w^s_h \in V^s_h , 
\end{array}
\end{equation}
for $s\in\{1,\ldots,m\}$, that results in the following linear system of equations:
\begin{equation} \label{E:lin-sys-full-SPD}
\left\{
\begin{array}{r l l}
\bold{A} \overline{u} + \bold{B}^T \overline{p} & =\overline{ \mathrm{f}}, \\[2pt]
\bold{B} \overline{u} - [\bold{\Sigma}_\varepsilon  \boldsymbol{\mathcal{B}}_{\mathcal{D}} + \bold{Q}] \overline{p}  & =\overline{0},
\end{array}
\right.\quad \overline{u} \in \mathbb{R}^N, \quad  \overline{p} \in \mathbb{R}^n,
\end{equation}
or equivalently,
\begin{equation} \label{E:lin-sys-z}
 \boldsymbol{\mathcal{A}}_\varepsilon \bold{z}_\varepsilon = \overline{ \mathrm{F}} ,
\end{equation}
with the {\it saddle point} matrix
\begin{equation} \label{E:lin-sys-full-matrix-1}
 \boldsymbol{\mathcal{A}}_\varepsilon = \begin{bmatrix} \bold{A} & \bold{B}^T \\ \bold{B} & - \bold{\Sigma}_\varepsilon  \boldsymbol{\mathcal{B}}_{\mathcal{D}} - \bold{Q}  \end{bmatrix}  \in \mathbb{R}^{(N+n)\times(N+n)},  \notag
\end{equation}
and vectors 
\[
\bold{z}_\varepsilon = 
\begin{bmatrix}
\overline{u} \\
\overline{p} 
\end{bmatrix} \in \mathbb{R}^{N+n},
\quad \overline{ \mathrm{F}}= \begin{bmatrix} \overline{ \mathrm{f}} \\ \overline{0} \end{bmatrix} \in \mathbb{R}^{N+n}.
\]

To provide the comprehensive description of the linear system   \eqref{E:lin-sys-full-SPD} or \eqref{E:lin-sys-z}, we introduce the following notations for the number of degrees of freedom in different parts of  $\Omega_h$. Let $N$ be the total number of nodes in $\Omega_h$, and 
$n$ be the number of nodes in $\overline{\mathcal{D}}_h$ so that 
\[
n=\sum_{s=1}^m n_s,
\]
where $n_s$ denotes the number of nodes in  $\overline{\mathcal{D}}^s_h$, and, finally, $n_{0}$ is the number of nodes in $\Omega_h \setminus \overline{\mathcal{D}_h}$, so that we have
\[
N=n_{0} + n .
\]
Then in \eqref{E:lin-sys-full-SPD}, the vector $\overline{u} \in \mathbb{R}^{N}$ has entries $u_i=u_h(x_i) $ with $x_i \in \Omega_h$. We count the entries of $\overline{u}$ in such a way that its first $n$ entries correspond to the nodes of $\overline{\mathcal{D}}_h$, and the remaining $n_0$ entries correspond to the nodes of $\Omega_h \setminus \overline{ \mathcal{D}}_h$. Entries of the first group can be further partitioned into $m$ subgroups such that there are $n_s$ entries in the $s^{\text{th}}$ group that corresponds to $\overline{\mathcal{D}}^s_h$, $s\in\{1,\ldots,m\}$.
Similarly, the vector $\overline{p}\in \mathbb{R}^{n}$ has entries $p_i=p_h(x_i) $ where $x_i \in \overline{ \mathcal{D}}_h$. Then we can write
\[
\mathbb{R}^{n}\ni\overline{p} =  \begin{bmatrix}
\overline{p}_1 \\
\vdots\\
\overline{p}_n
\end{bmatrix} , \quad \mbox{where }~ \overline{p}_s \in \mathbb{R}^{n_s}, \quad s\in\{1,\ldots,m\}.
\]

The symmetric positive definite matrix $\bold{A} \in\mathbb{R}^{N\times N}$ of \eqref{E:lin-sys-full-SPD}  is the stiffness matrix that arises from the discretization of the Laplace operator with the homogeneous Dirichlet boundary conditions on $\Gamma$, that is,
\begin{equation} \label{E:A-def}
(\bold{A} \overline{u} , \overline{v} ) = \int_{\Omega_h} \nabla u_h \cdot \nabla v_h ~dx, \quad\mbox{where}\quad \overline{u} , \overline{v} \in \mathbb{R}^N, \quad
u_h , v_h \in V_h,
\end{equation}
where $(\cdot,\cdot)$ is the standard dot-product of vectors.
With the above orderings, the matrix $\bold{A} $ of \eqref{E:A-def} can be presented as $2\times 2$ block-matrix 
\begin{equation} \label{E:matr-A}
\bold{A} = 
\begin{bmatrix}
\mathrm{A}_{\mathcal{D}\mathcal{D}} & \mathrm{A}_{\mathcal{D}0} \\
\mathrm{A}_{0\mathcal{D}} & \mathrm{A}_{00} 
\end{bmatrix},
\end{equation}
where the block $\mathrm{A}_{\mathcal{D}\mathcal{D}} \in  \mathbb{R}^{n\times n}$ corresponds to the inclusions $\overline{\mathcal{D}}^s_h$, $s\in \{1,\ldots,m\}$, the  block $\mathrm{A}_{00}\in  \mathbb{R}^{n_0\times n_0}$ corresponds to the region outside of $\overline{\mathcal{D}}_h$, and the entries of $\mathrm{A}_{\mathcal{D}0} \in  \mathbb{R}^{n\times n_0}$ and $\mathrm{A}_{0\mathcal{D}}=\mathrm{A}_{\mathcal{D}0}^T$ are assembled from entries associated with both $\overline{\mathcal{D}}_h$ and $\Omega_h \setminus \overline{\mathcal{D}}_h$. 

The matrix $\boldsymbol{\mathcal{B}}_{\mathcal{D}} \in \mathbb{R}^{n\times n}$ in \eqref{E:lin-sys-full-SPD}, that corresponds to the highly conducting inclusions, is the $m\times m$ block-diagonal matrix
\begin{equation}\label{E:matrix-B_D}
\boldsymbol{\mathcal{B}}_{\mathcal{D}} = 
\text{diag}~ (\mathrm{B_1}, \ldots, \mathrm{B_m}),
\end{equation}
whose blocks $\mathrm{B_s}\in \mathbb{R}^{n_s \times n_s}$ are  defined  by
\begin{equation} \label{E:B-def}
(\mathrm{B_s} \overline{u} , \overline{v} ) = \int_{\mathcal{D}^s} \nabla u_h \cdot \nabla v_h ~dx, \quad\mbox{where}\quad \overline{u} , \overline{v} \in \mathbb{R}^{n_s}, \quad u_h , v_h \in V^s_h.
\end{equation}
Note that the matrix $\mathrm{B_s}$ is the stiffness matrix in the discretization of the Laplace operator in the domain $\mathcal{D}^s$ with the Neumann boundary conditions on $\Gamma_s$, $s\in\{1,\ldots,m\}$.  Also, remark that each matrix $\mathrm{B_s}$ is {\it positive semidefinite} with 
\begin{equation} \label{E:ker-Bi}
\ker \mathrm{B_s} = \mbox{span}
\left\{\overline{e}_s \right\}, \quad \mbox{where} \quad \overline{e}_s=\begin{bmatrix}
1 \\
\vdots\\
1
\end{bmatrix} \in \mathbb{R}^{n_s}.
\end{equation}
To this end, 
\[
\dim \ker ~ \boldsymbol{\mathcal{B}}_{\mathcal{D}}  = m.
\]
Then, the matrix $\bold{B} \in \mathbb{R}^{n\times N}$ of \eqref{E:lin-sys-full-SPD} is written in the block form  as
\begin{equation} \label{E:matr-B}
\bold{B} = \begin{bmatrix} \boldsymbol{\mathcal{B}}_{\mathcal{D}}  & \bold{0}\end{bmatrix},
\end{equation}
with zero-matrix $\bold{0} \in \mathbb{R}^{n\times n_0}$ and $\boldsymbol{\mathcal{B}}_{\mathcal{D}} \in \mathbb{R}^{n\times n}$. The vector $\overline{ \mathrm{f}}  \in \mathbb{R}^{N}$ of  \eqref{E:lin-sys-full-SPD}  is defined in a similar way by
\[
(\overline{ \mathrm{f}},\overline{v}) = \int_{\Omega_h} f v_h ~dx, \quad\mbox{where}\quad \overline{v} \in \mathbb{R}^N, \quad v_h \in V_h.
\]
With all that, the first equation of \eqref{E:FEM-form} results in the first equation of \eqref{E:lin-sys-full-SPD}.
Now denote
\begin{equation*} 
\bold{\Sigma}_\varepsilon =  
\text{diag}~ (\varepsilon_1\mathrm{I_1}, \ldots, \varepsilon_m\mathrm{I_m}),
\end{equation*}
with $\mathrm{I_s} \in \mathbb{R}^{n_s\times n_s}$ being the identity matrix. Finally, we construct the matrix $\bold{Q}$  in \eqref{E:FEM-form} using
\begin{equation} \label{E:Q}
\bold{Q} =\text{diag}~ ( \mathrm{Q_1}, \ldots,  \mathrm{Q_m}), 
\end{equation}
whose blocks $\mathrm{Q_s}  \in \mathbb{R}^{n_s\times n_s}$, $s\in\{1,\ldots,m\}$, are defined by
\begin{equation} \label{E:Q-def}
(\mathrm{Q_s} \overline{p} , \overline{q} ) =  \frac{1}{|\mathcal{D}^s_h|}   \left[ \int_{\mathcal{D}^s_h} p_h~dx \right]  \left[ \int_{\mathcal{D}^s_h} q_h~dx \right], ~~\mbox{where}~ \overline{p} , \overline{q} \in \mathbb{R}^{n_s}, ~~ p_h, ~ q_h \in V^s_h.
\end{equation}
As would be evident from below considerations, another way of writing the matrix $\mathrm{Q_s}$ is via
\begin{equation} \label{E:evector-1}
\mathrm{Q_s} =  \frac{1}{d_s^2}  \left[\mathrm{M_s} \overline{w}^1_s \otimes \mathrm{M_s} \overline{w}^1_s  \right], \quad \mbox{where} \quad d_s=\left|\mathcal{D}^s_h \right|^{1/2},
\quad \mbox{and} \quad 
\overline{w}^1_s := \frac{1}{d_s} \overline{e}_s   \in \mathbb{R}^{n_s},  
\end{equation}
and $\mathrm{M_s} \in \mathbb{R}^{n_s \times n_s}$ is the {\it mass matrix} associated with the inclusion $\mathcal{D}^s$ and given by 
\begin{equation} \label{E:mass1}
 (\mathrm{M_s} \overline{p}_s , \overline{q}_s ) = \int_{\mathcal{D}^s_h} p^s_h \, q^s_h ~dx, \quad \mbox{for all } \ \overline{p}_s,\overline{q}_s, \in \mathbb{R}^{n_s}, \quad  p^s_h,~ q^s_h \in V^s_h, \ s\in \{1,\ldots,m\}.
\end{equation}
In \eqref{E:evector-1}, $\overline{p} \otimes \overline{q} = \overline{p} \, \overline{q}^T$ denotes the {\it outer product} of vectors $ \overline{p} $ and $\overline{q}$. 
The matrix $\mathrm{Q_s}$ is a symmetric and positive semidefinite rank-one matrix generated by the $\mathrm{M_s}$-normal vector $\overline{w}^1_s $, that is, $( \mathrm{M_s} \overline{w}^1_s , \overline{w}^1_s) = 1$, $s\in \{1,\ldots,m\}$.

With \eqref{E:B-def}--\eqref{E:mass1},  the second equation of \eqref{E:var-form-ours} yields
the second equation in the system  \eqref{E:lin-sys-full-SPD}. Note that with \eqref{E:matr-A}, the symmetric and indefinite matrix 
$ \boldsymbol{\mathcal{A}}_\varepsilon$ defined in \eqref{E:lin-sys-z} is then
\begin{equation} \label{E:lin-sys-full-matrix-2}
\boldsymbol{\mathcal{A}}_\varepsilon 
= \begin{bmatrix} \mathrm{A}_{\mathcal{D}\mathcal{D}} & \mathrm{A}_{\mathcal{D}0} & \boldsymbol{\mathcal{B}}_{\mathcal{D}} \\ \mathrm{A}_{0\mathcal{D}} & \mathrm{A}_{00} &  \bold{0}^T \\ \boldsymbol{\mathcal{B}}_{\mathcal{D}}  & \bold{0} & - \bold{\Sigma}_\varepsilon  \boldsymbol{\mathcal{B}}_{\mathcal{D}} - \bold{Q}  \end{bmatrix}.
\end{equation}
This concludes the derivation of the saddle point formulation \eqref{E:lin-sys-full-SPD}. 
Clearly, there exists a unique solution $\overline{u} \in \mathbb{R}^N$, $\overline{p} \in \mathbb{R}^n$, or equivalently, $\bold{z}_\varepsilon \in \mathbb{R}^{N+n}$.

System \eqref{E:lin-sys-full-SPD} was proposed in \cite{gkk18,kuz00,kuz09}  for the case when $\bold{Q}=\bold{0}$, where it was also demonstrated that  \eqref{E:lin-sys-full-SPD} can be derived in a purely algebraic way.

\section{Preconditioned Iterative Methods} \label{S:methods}

In this paper, we consider and investigate three iterative methods for solving system \eqref{E:lin-sys-z}. 
The first one is the preconditioned conjugate gradient method or preconditioned Uzawa (PU) for the Schur complement system
\begin{equation}   \label{E:Uzawa-syst}
\bold{S}_\varepsilon  \, \overline{p} 
=\overline{ \mathrm{g}}  = : \bold{B}\bold{A}^{-1} \, \overline{ \mathrm{f}},
\end{equation}
where 
\begin{equation} \label{E:schur}
\bold{S}_\varepsilon : = \bold{\Sigma}_\varepsilon  \boldsymbol{\mathcal{B}}_{\mathcal{D}} + \bold{Q} + \bold{B}\bold{A}^{-1}\bold{B}^T,
\end{equation} 
with the preconditioner 
\begin{equation}   \label{E:Uzawa-prec}
\mathcal{H}_{\mathrm{S}} =  [\boldsymbol{\mathcal{B}}_{\mathcal{D}} + \bold{Q}]^{-1} \in \mathbb{R}^{n \times n} .
\end{equation}

The second method is the preconditioned Lanzcos (PL) method with the preconditioner
\begin{equation}   \label{E:block-prec}
\boldsymbol{\mathcal{H}} = 
\begin{bmatrix}
\mathcal{H}_{\mathrm{A}} & 0\\
0 & \mathcal{H}_{\mathrm{S}}
\end{bmatrix}, 
\end{equation} 
where $\mathcal{H}_{\mathrm{A}} \in \mathbb{R}^{N \times N}$ is a given symmetric positive definite matrix introduced below, and $\mathcal{H}_{\mathrm{S}}$ is the same as in \eqref{E:Uzawa-prec}. 

The third method is the preconditioned conjugate gradient (PCG) method with the preconditioner $\boldsymbol{\mathcal{H}}$ defined in (\ref{E:block-prec}) for a modified system obtained from \eqref{E:lin-sys-z} as follows:
\begin{equation}\label{E:mod-syst_PCG}
\bold{K}_\varepsilon 
{\bold{z}}_\varepsilon = \boldsymbol{\mathcal{G}}_\varepsilon ,
\end{equation}
where
\begin{equation}\label{E:mod-syst_PCG_k_e}
\bold{K}_\varepsilon = \boldsymbol{\mathcal{A}}_\varepsilon \boldsymbol{\mathcal{H}}   \boldsymbol{\mathcal{A}}_\varepsilon, \quad 
 \boldsymbol{\mathcal{G}}_\varepsilon= \boldsymbol{\mathcal{A}}_\varepsilon  \boldsymbol{\mathcal{H}} \, \overline{ \mathrm{F}}.
\end{equation}

\subsection{ Preconditioned Uzawa Method} \label{S:Uzawa}

The preconditioned Uzawa algorithm combined with the PCG method is well known, see e.g. \cite{bpx90,rg84}. It is defined by 
\begin{equation}\label{E:PCG_Uzawa-1}
\overline{p}^k = \overline{p}^{k-1} - \beta_k \overline{\xi}_k,  \quad k=1,2,\ldots,
\end{equation}
where
\begin{equation}\label{E:PCG_Uzawa-0}
\overline{\xi}_k = \begin{cases} 
\mathcal{H}_{\mathrm{S}} (\bold{S}_\varepsilon  \overline{p}^{0} - \overline{ \mathrm{g}} ), &k=1 \\
\mathcal{H}_{\mathrm{S}} (\bold{S}_\varepsilon  \overline{p}^{k-1} - \overline{ \mathrm{g}} ) - \alpha_k  \overline{\xi}_{k-1}, &k \geq 2,
\end{cases}
\end{equation}
and
\begin{equation}\label{E:PCG_Uzawa-2}
\beta_k = \frac{ \left(  \bold{S}_\varepsilon  \overline{p}^{k-1} -  \overline{ \mathrm{g}} , \overline{\xi}_{k} \right) }{(\bold{S}_\varepsilon  \overline{\xi}_{k},  \overline{\xi}_{k})}, \quad
\alpha_k = \frac{ \left( \mathcal{H}_{\mathrm{S}} (\bold{S}_\varepsilon  \overline{p}^{k-1} -  \overline{ \mathrm{g}} ) , \bold{S}_\varepsilon   \overline{\xi}_{k-1}\right) }{(\bold{S}_\varepsilon \overline{\xi}_{k-1},  \overline{\xi}_{k-1})}, \quad k=1,2,\ldots
\end{equation}
Here $\overline{p}^0$ is an initial guess, and $\mathcal{H}_{\mathrm{S}}$ is given by (\ref{E:Uzawa-prec}).

Denote  by $\overline{p}^*$ the solution of \eqref{E:Uzawa-syst}, 
then the convergence estimate for \eqref{E:mod-syst_PCG}-\eqref{E:PCG_Uzawa-1} is given by (see  \cite{axel,km74}):
\begin{equation*}
\| \overline{p}^k - \overline{p}^* \|_{\bold{S}_\varepsilon} 
\leqslant 
\frac{1}{C_k \left(\frac{b + a}{b - a}\right)} \| \overline{p}^0 - \overline{p}^* \|_{\bold{S}_\varepsilon}, \quad k= 0,1,2,\ldots,
\end{equation*}
where $\| \cdot \|_{\bold{S}_\varepsilon}$ is the elliptic norm generated by the matrix $\bold{S}_\varepsilon$, $C_k\left(t\right)$ is the Chebyshev polynomial of degree $k$, and $b$ and $a$ are the estimates from above and from below for the eigenvalues of the matrix $\mathcal{H}_{\mathrm{S}} \bold{S}_\varepsilon$, respectively.

To investigate the eigenvalue problem
\begin{equation*}
\mathcal{H}_{\mathrm{S}} \bold{S}_\varepsilon \overline{\psi} = \mu \,  \overline{\psi},
\end{equation*}
we observe that
\begin{equation}\label{E:eigenval_Q1}
\mathcal{H}_{\mathrm{S}} \boldsymbol{\mathcal{B}}_{\mathcal{D}} = \bold{I} - \tilde{\bold{Q}}, 
\end{equation}
\begin{equation}\label{E:eigenval_Q2}
\mathcal{H}_{\mathrm{S}} \bold{Q} = \tilde{\bold{Q}},
\end{equation}
where $ \bold{\tilde{Q}}$ is $m \times m$ block diagonal matrix:
\begin{equation*}
\tilde{ \bold{Q}}  =\text{diag}~ ( \tilde{\mathrm{Q}}_1, \ldots,  \tilde{\mathrm{Q}}_m),
\end{equation*}             
with $\mathrm{M_s}$-orthogonal projectors
\begin{equation*}
\tilde{\mathrm{Q}}_s = \overline{w}^1_s \otimes \left(\mathrm{M_s} \overline{w}^1_s \right) \in \mathbb{R}^{n_s \times n_s}, \quad s\in \{1,\ldots,m\},
\end{equation*}                                  
where $\overline{w}^1_s$  and $\mathrm{M}_s$ were introduced in \eqref{E:evector-1} and \eqref{E:mass1}, respectively.

\begin{remark} \label{R:rk1}
It follows from (\ref{E:eigenval_Q1}), (\ref{E:eigenval_Q2}) that implementation of the matrix-vector products  $\mathcal{H}_{\mathrm{S}} \boldsymbol{\mathcal{B}}_{\mathcal{D}} \, \overline{y}$ and $\mathcal{H}_{\mathrm{S}} \bold{Q}\, \overline{y}$ requires only $2n$ arithmetical operations for any vector $\overline{y} \in  \mathbb{R}^n$, that is, we do not need to solve a system with the matrix $\boldsymbol{\mathcal{B}}_{\mathcal{D}} + \bold{Q}$.
\end{remark}

Simple algebraic analysis, see e.g. \cite{gkk18,kuz09}, shows that
\begin{equation}\label{E:a_0_estim}
a \geqslant \min \{ a_{0}+ \varepsilon_{\min} ;1\}
\end{equation}
and 
\begin{equation}\label{E:b_0_estim}
b \geqslant \max \{ b_{0}+\varepsilon_{\max} ;1\}
\end{equation}
where $a_{0}>0$ and $b_{0}$ are estimates from below and above, respectively, for the eigenvalues of the matrix
\begin{equation*} 
\mathcal{H}_{\mathrm{S}} \boldsymbol{\mathrm{S}}_0 \equiv \mathcal{H}_{\mathrm{S}} \bold{B} \bold{A}^{-1} \bold{B}^{T}  
\end{equation*}
that is, $\boldsymbol{\mathrm{S}}_0 =\boldsymbol{\mathrm{S}}_\varepsilon$ when $\varepsilon_1= \ldots\varepsilon_m=0$, and $\varepsilon_{\min} = \min \limits_{1 \leqslant t \leqslant m} \varepsilon_t$, $\varepsilon_{\max} = \max \limits_{1 \leqslant t \leqslant m} \varepsilon_t$. 
The values of  $a_{0}$ and $b_{0}$ will be derived in Section \ref{S:eig_est}.

\subsection{Preconditioned Lanczos Method}  \label{S:PL}

Preconditioned Lanczos  method for systems with symmetric indefinite matrices was proposed in late 1960s, see \cite{km74} and references therein. 
In this paper, we consider PL method for the saddle-point system (\ref{E:lin-sys-z}) preconditioned by a symmetric positive definite matrix $\boldsymbol{\mathcal{H}}$ of \eqref{E:block-prec} with some 
given symmetric positive definite matrix $\mathcal{H}_{\mathrm{A}}$ introduced below, and $\mathcal{H}_{\mathrm{S}}$ defined by \eqref{E:Uzawa-prec}. The PL method is as follows, see e.g. \cite{km74}:
\begin{equation}
\overline{z}^k = \overline{z}^{k-1} - \beta_k \overline{\xi}_k,  \quad  k= 1,2,\ldots,
\end{equation}
where
\begin{equation}
\label{E:Prec_lanc-0}
\overline{\xi}_k = \begin{cases} 
\boldsymbol{\mathcal{H}}(\boldsymbol{\mathcal{A}}_{\varepsilon}\overline{z}^0 -  \overline{ \mathrm{F}}  ), &k=1 \\
\boldsymbol{\mathcal{H}}\boldsymbol{\mathcal{A}}_{\varepsilon}\overline{\xi}_1 - \alpha_2 \overline{\xi}_1, & k=2\\
\boldsymbol{\mathcal{H}}\boldsymbol{\mathcal{A}}_{\varepsilon}\overline{\xi}_{k-1} - \alpha_k \overline{\xi}_{k-1} - \gamma_k \overline{\xi}_{k-2}, & k\geq 3,
\end{cases}
\end{equation}
and
\begin{equation}
\alpha_k = \frac{(\boldsymbol{\mathcal{A}}_{\varepsilon} \boldsymbol{\mathcal{H}} \boldsymbol{\mathcal{A}}_{\varepsilon} \overline{\xi}_{k-1}, \boldsymbol{\mathcal{H}} \boldsymbol{\mathcal{A}}_{\varepsilon} \overline{\xi}_{k-1})}{(\boldsymbol{\mathcal{A}}_{\varepsilon} \overline{\xi}_{k-1},  \boldsymbol{\mathcal{H}} \boldsymbol{\mathcal{A}}_{\varepsilon} \overline{\xi}_{k-1})}, \qquad  
\gamma_k = \frac{(\boldsymbol{\mathcal{A}}_{\varepsilon} \boldsymbol{\mathcal{H}}\boldsymbol{\mathcal{A}}_{\varepsilon} \overline{\xi}_{k-1},  \boldsymbol{\mathcal{H}}\boldsymbol{\mathcal{A}}_{\varepsilon} \overline{\xi}_{k-2})}{(\boldsymbol{\mathcal{A}}_{\varepsilon} \overline{\xi}_{k-2},  \boldsymbol{\mathcal{H}} \boldsymbol{\mathcal{A}}_{\varepsilon} \overline{\xi}_{k-2})}, \quad  k= 1,2,\ldots,
\end{equation}
\begin{equation}
\beta_k = \frac{( \boldsymbol{\mathcal{A}}_{\varepsilon}\overline{z}^{k-1} -  \overline{ \mathrm{F}} , \boldsymbol{\mathcal{H}} \boldsymbol{\mathcal{A}}_{\varepsilon} \overline{\xi}_{k})}{(\boldsymbol{\mathcal{A}}_{\varepsilon} \overline{\xi}_{k},  \boldsymbol{\mathcal{H}} \boldsymbol{\mathcal{A}}_{\varepsilon} \overline{\xi}_{k})},  \quad k= 1,2,\ldots
\label{E:Prec_lanc-3}
\end{equation}
Let $\overline{z}^0$ be an initial guess, $\overline{z}^*$ the solution of \eqref{E:lin-sys-z}, then, the following convergence estimate holds
\begin{equation}\label{E:eigenval_est_PL}
\| \overline{z}^k - \overline{z}^* \|_{\bold{K}_\varepsilon} 
\leqslant 
\frac{1}{C_{k/2} \left(\frac{b^2 + a^2}{b^2 - a^2}\right)} \| \overline{z}^0 - \overline{z}^* \|_{\bold{K}_\varepsilon}, \quad k=2,4,\ldots
\end{equation}
see \cite{km74},
where $\bold{K}_\varepsilon$ is given by \eqref{E:mod-syst_PCG_k_e} and  $\| \cdot \|_{\bold{K}_\varepsilon}$ is the elliptic norm generated by the matrix $\bold{K}_\varepsilon =\bold{K}_\varepsilon^T > 0 $. Here $C_{k/2}$ is the Chebyshev polynomial of degree $k/2$, and $b^2$ and $a^2>0$ are estimates from above and from below for eigenvalues of the matrix  $\left(\boldsymbol{\mathcal{H}}   \boldsymbol{\mathcal{A}}_\varepsilon\right)^2 $, respectively.

\subsection{Preconditioned Conjugate Gradient method}    \label{S:PCG}

The preconditioned conjugate gradient method with the preconditioner $\boldsymbol{\mathcal{H}}$ defined by \eqref{E:block-prec},  we apply to system (\ref{E:mod-syst_PCG}):
\begin{equation}
\overline{z}^k = \overline{z}^{k-1} - \beta_k \overline{\xi}_k,  \quad k=1,2,\ldots ,
\end{equation}
where
\begin{equation}
\overline{\xi}_k = \begin{cases} 
\boldsymbol{\mathcal{H}}  (\bold {K}_\varepsilon \overline{z}^{0} - \boldsymbol{\mathcal{G}}_\varepsilon ), &k=1 \\
\boldsymbol{\mathcal{H}}  (\bold{K}_\varepsilon \overline{z}^{k-1} - \boldsymbol{\mathcal{G}}_\varepsilon ) - \alpha_k  \overline{\xi}^{k-1}, &k \geq 2,
\end{cases}
\label{E:prec_pcg-0}
\end{equation}
and
\begin{equation}
\alpha_k = \frac{( \boldsymbol{\mathcal{H}}[ \bold{K}_\varepsilon \overline{z}^{k-1} -  \boldsymbol{\mathcal{G}}_\varepsilon ],  \bold{K}_\varepsilon \overline{\xi}_{k-1})}{(  \bold{K}_\varepsilon \overline{\xi}_{k-1}, \overline{\xi}_{k-1})}
, \qquad  
\beta_k = \frac{( \bold{K}_\varepsilon \overline{z}^{k-1} -  \boldsymbol{\mathcal{G}}_\varepsilon, \overline{\xi}_{k})}{( \bold{K}_\varepsilon \overline{\xi}_{k}, \overline{\xi}_{k})}, \quad k=1,2,\ldots 
\label{E:prec_pcg-2}
\end{equation}
The convergence estimate for the method is as follows, see \cite{axel,km74}: 
\begin{equation}
\| \overline{z}^k - \overline{z}^* \|_{\bold{K}_\varepsilon} 
\leqslant 
\frac{1}{C^{k} \left(\frac{b^2 + a^2}{b^2 - a^2}\right)} \| \overline{z}^0-\overline{z}^* \|_{\bold{K}_\varepsilon}, \ k=1,2,\ldots
\end{equation}
with the same matrix  $\bold{K}_\varepsilon$ defined in \eqref{E:mod-syst_PCG_k_e} and values $a^2$ and $b^2$ as in \eqref{E:eigenval_est_PL}.

\section{Eigenvalue estimates}\label{S:eig_est}
\subsection{Eigenvalue estimates for the matrix $\mathcal{H}_{\mathrm{S}} \boldsymbol{\mathrm{S}}_0 $} \label{S:eig_est_hs}

Consider the eigenvalue problem 

\begin{equation}\label{E:eigv_prob_init}
 \boldsymbol{\mathrm{S}}_0 \overline{\psi}_{\mathcal{D}} = \mu \, \mathcal{H}_{\mathrm{S}}^{-1} \overline{\psi}_{\mathcal{D}},
\end{equation}
where
\begin{equation}
 \boldsymbol{\mathrm{S}}_0 = \bold{B} \bold{A}^{-1} \bold{B}^T = \boldsymbol{\mathcal{B}}_{\mathcal{D}} \mathrm{S}^{-1}_{00} \boldsymbol{\mathcal{B}}_{\mathcal{D}},
\end{equation}
\begin{equation}
 \mathcal{H}_{\mathrm{S}}^{-1}  =  \boldsymbol{\mathcal{B}}_{\mathcal{D}} +  \bold {Q},
\end{equation}
and
\begin{equation*}
\mathrm{S}_{00}= \mathrm{A}_{\mathcal{D}\mathcal{D}}  - \mathrm{A}_{\mathcal{D}0} \mathrm{A}_{00}^{-1}\mathrm{A}_{0\mathcal{D}}	
\end{equation*}
is the Schur complement of $\mathrm{A}_{00}$. It is obvious that $\lambda = 1$ if $\overline{\psi}_{\mathcal{D}} \in \ker \boldsymbol{\mathcal{B}}_\mathcal{D}$, 
and $\bold {Q} \overline{\psi}_{\mathcal{D}}=\boldsymbol{0}$ for any $\lambda \neq 1$, that is,  $\overline{\psi}_{\mathcal{D}}$ is $\bold{M}$-orthogonal to  $\ker \boldsymbol{\mathcal{B}}_{\mathcal{D}}$, with
\[
\bold{M} = \mbox{diag} \left \{ \mathrm{M_1}, \ldots, \mathrm{M_m} \right\}, \quad\mbox{where } \mathrm{M_s} \mbox{ is given by } \eqref{E:mass1}, \quad s\in \{1,\ldots,m\}. 
\]
Thus, to derive $a_0$ and $b_0$ of (\ref{E:a_0_estim})-(\ref{E:b_0_estim}),  instead of \eqref{E:eigv_prob_init},  we can consider the following eigenvalue problem
\begin{equation}\label{E:eig_val_main_eq}
\boldsymbol{\mathrm{S}}_0 \overline{\psi}_{\mathcal{D}} =\mu \, \boldsymbol{\mathcal{B}}_{\mathcal{D}} \overline{\psi}_{\mathcal{D}},
\end{equation}
under the condition $ \left( \bold{M} \overline{\psi}_{\mathcal{D}} ,\overline{w}\right)=0$ for all $\overline{w} \in \ker  \boldsymbol{\mathcal{B}}_{\mathcal{D}} $.
 
Let $\mu$ be an eigenvalue of (\ref{E:eig_val_main_eq}) and $\overline{\psi}_{\mathcal{D}}$ a corresponding eigenvector. Then,
\begin{equation} \label{E:raylegh}
\mu = \frac{\left(\boldsymbol{\mathrm{S}}_0 \overline{\psi}_{\mathcal{D}},\overline{\psi}_{\mathcal{D}} \right)}{\left( \boldsymbol{\mathcal{B}}_{\mathcal{D}}  \overline{\psi}_{\mathcal{D}},\overline{\psi}_{\mathcal{D}} \right)} = \frac{\left(\bold{A} \overline{\psi},\overline{\psi} \right)}{\left( \boldsymbol{\mathcal{B}}_{\mathcal{D}}  \overline{\psi}_{\mathcal{D}},\overline{\psi}_{\mathcal{D}} \right)} = \frac{ \int \limits_{\Omega_h} \left| \nabla  \psi_h \right|^2 dx}{ \int \limits_{\mathcal{D}_h} \left| \nabla  \psi_h \right|^2  dx},
\end{equation}
where 
\begin{equation*}
\overline{\psi} = \begin{bmatrix}
\overline{\psi}_{\mathcal{D}} \\
\overline{\psi}_0
\end{bmatrix}, \quad
\mbox{such that} \quad \mathrm{A}_{0\mathcal{D}} \overline{\psi}_{\mathcal{D}} + \mathrm{A}_{00}\overline{\psi}_0=0, 
\end{equation*}
and ${\psi}_h \in V_h$. The vector $\overline{\psi}_0 \in \mathbb{R}^{n_0}$ corresponds to a FEM function ${\psi}_{0,h} \in V_h|_{\Omega_h\setminus\mathcal{D}_h}$
called the {\it continuous $h$-harmonic extension} of ${\psi}_{\mathcal{D},h} \in V_h|_{\mathcal{D}_h}$ from $\mathcal{D}_h$ into $\Omega_h \setminus \mathcal{D}_h$, where the FEM function ${\psi}_{\mathcal{D},h}$ corresponds to the vector $\overline{\psi}_{\mathcal{D}}  \in \mathbb{R}^n$. Note that ${\psi}_{0,h}$ is the solution of the following variational finite element problem:
\begin{equation*}
\begin{array}{l l l}
& \displaystyle  \mbox{  Find  } & u_h \in V_h|_{\Omega_h\setminus\mathcal{D}_h} ~ \mbox{ satisfying } ~ u_h =  {\psi}_{\mathcal{D},h}  ~ \mbox{ on } ~ \partial \mathcal{D}_h ~ \mbox{ such that  } \\[2pt]
& & \displaystyle   \int \limits_{\Omega_h \setminus \mathcal{D}_h} \left| \nabla u_h \right|^2 dx = \min_{ \substack{v_h \in V_h|_{\Omega_h\setminus\mathcal{D}_h}\\ v_h|_{\partial \mathcal{D}_h } =  {\psi}_{\mathcal{D},h} }} ~\int \limits_{\Omega_h\setminus \mathcal{D}_h} \left| \nabla v_h \right|^2 dx.
\end{array}
\end{equation*}
From now on, we will write $\mathcal{D}$ instead of $\mathcal{D}_h$, and $\mathcal{D}^s$ instead of $\mathcal{D}_h^s$, $s\in \{1,\ldots,m\}$, since they are the same due to the above assumptions.

To estimate the value of $a_0$ in (\ref{E:a_0_estim}) from below we consider the eigenvalue problem \eqref{E:eig_val_main_eq} using the spectral decomposition of $\mathrm{B_s} \in \mathbb{R}^{n_s}$, $s\in \{1,\ldots,m\}$ that comes from 
\begin{equation}
\mathrm{B_s} \ol{w} = \lambda \mathrm{M_s} \ol{w},
\label{e:kuz_eigen_mass}
\end{equation}
that is,
\begin{equation}
\mathrm{B_s} = \mathrm{M_s} \mathrm{W_s} \Lambda_s \mathrm{W_s}^T \mathrm{M_s},
\label{e:kuz_spectral_decomp}
\end{equation}
with
\begin{equation*}
\mathrm{W_s} = \left[\ol{w}_s^1, \ldots,\ol{w}_s^{n_s} \right], \quad  \mbox{and} \quad 
\Lambda_s = \mbox{diag} \left\{ \lambda_s^1, \ldots, \lambda_s^{n_s} \right\},
\end{equation*}
where $0=\lambda_s^1 < \lambda_s^2 \leqslant \ldots \leqslant \lambda_s^{n_s}$ are the eigenvalues in \eqref{e:kuz_eigen_mass} and $\ol{w}_s^1, \ldots,\ol{w}_s^{n_s}$ are the corresponding $\mathrm{M_s}$-orthonormal eigenvectors, $s\in \{1,\ldots,m\}$. 
We define the matrices
\begin{equation}
\mathrm{\hat{B}_s} = \mathrm{M_s}^{\frac{1}{2}} \mathrm{W_s} \Lambda_s \mathrm{W_s}^T  \mathrm{M_s}^{\frac{1}{2}},
\label{e:kuz_B_matrix_decomp}
\end{equation}
and
\begin{equation}
\hat{\mathrm{B}}_s^{\frac{1}{2}} = \mathrm{M_s}^{\frac{1}{2}} \mathrm{W_s} \Lambda_s^{\frac{1}{2}} \mathrm{W_s}^T  \mathrm{M_s}^{\frac{1}{2}},
\label{e:kuz_B_matrix_decomp_half}
\end{equation}
It is obvious that $\hat{\mathrm{B}}_s^{\frac{1}{2}} $ are symmetric positive semidefinite matrices and $\hat{\mathrm{B}}_s^{\frac{1}{2}} \hat{\mathrm{B}}_s^{\frac{1}{2}}  =\hat{\mathrm{B}}_s$, $s\in \{1,\ldots,m\}$. 
Also note that $\ol{w}_s^1 \in \ker \mathrm{B}_s$ and is precisely the one that is given by \eqref{E:evector-1}. 
In addition, we define the matrices 
\begin{equation*}
\hat{\mathrm{B}}_{d, s}^{\frac{1}{2}} = \hat{\mathrm{B}}_s^{\frac{1}{2}} + \frac{1}{ d_s} \mathrm{M}_s^{\frac{1}{2}} \ol{w}_s^1 \otimes  \mathrm{M}_s^{\frac{1}{2}}\ol{w}_s^1,
\end{equation*}
where $d_s$, $s\in \{1,\ldots,m\}$, was introduced in  \eqref{E:evector-1}. 
Straightforward multiplications show that 
\begin{equation*}
\hat{\mathrm{B}}_{s}^{\frac{1}{2}}  \hat{\mathrm{B}}_{d, s}^{\frac{1}{2}} = \hat{\mathrm{B}}_{d, s}^{\frac{1}{2}}  \hat{\mathrm{B}}_s^{\frac{1}{2}} = \hat{\mathrm{B}}_s, \quad s\in \{1,\ldots,m\}.
\end{equation*}
The latter observation shows that  eigenvalue problem \eqref{E:eig_val_main_eq} is equivalent to the eigenvalue problem
\begin{equation}
\bold{M}^{\frac{1}{2}} \hat{\boldsymbol {\mathcal{B}}}^{\frac{1}{2}}  \hat{\boldsymbol {\mathcal{B}}}_d^{\frac{1}{2}}  \bold{M}^{\frac{1}{2}} \mathrm{S}_{00}^{-1} \bold{M}^{\frac{1}{2}}  \hat{\boldsymbol {\mathcal{B}}}_d^{\frac{1}{2}} \hat{\boldsymbol {\mathcal{B}}}^{\frac{1}{2}} \bold{M}^{\frac{1}{2}} \ol{w} = \mu \, \boldsymbol{\mathcal{B}}_{\mathcal{D}} \ol{w},
\label{e:kuz_equival_eigenval_pr}
\end{equation}
and
\begin{equation*}
\hat{\boldsymbol {\mathcal{B}}}^{\frac{1}{2}} = \mbox{diag} \left (\hat{\mathrm{B}}_1^{\frac{1}{2}}, \ldots, \hat{\mathrm{B}}_m^{\frac{1}{2}} \right), \quad 
\hat{\boldsymbol {\mathcal{B}}}_{d}^{\frac{1}{2}} = \mbox{diag} \left ( \hat{\mathrm{B}}_{d,1}^{\frac{1}{2}}, \ldots, \hat{\mathrm{B}}_{d,m}^{\frac{1}{2}} \right), 
\end{equation*}
are $m \times m $ block diagonal matrices. It is easy to see that the minimal eigenvalue in (\ref{e:kuz_equival_eigenval_pr}) is bounded from below by the minimal eigenvalue of the matrix
\begin{equation*}
\hat{\boldsymbol {\mathcal{B}}}_d^{\frac{1}{2}} \bold{M}^{\frac{1}{2}}\mathrm{S}_{00}^{-1} \bold{M}^{\frac{1}{2}} \hat{\boldsymbol{\mathcal{B}}}_d^{\frac{1}{2}}
\end{equation*}
which is equal to the minimal eigenvalue of the similar matrix $\mathrm{S}_{00}^{-1} \boldsymbol{\mathcal{B}}_{d}$ with
\begin{equation*}
\boldsymbol{\mathcal{B}}_d = \bold{M}^{\frac{1}{2}} \hat{\boldsymbol {\mathcal{B}}}_d \bold{M}^{\frac{1}{2}} = \mbox{diag} \left( \mathrm{B}_1+\mathrm{Q_1}, \ldots, \mathrm{B}_s+ \mathrm{Q_m} \right) , 
\end{equation*}
where $\mathrm{Q_s}$, $s\in \{1,\ldots,m\}$ is defined in \eqref{E:evector-1}.
If $\left(\mu, \ol{w} \right) $ is an eigenpair of the matrix  $\mathrm{S}_{00}^{-1} \boldsymbol{\mathcal{B}}_{d}$, then  similar to \eqref{E:raylegh}, we obtain

\begin{equation}
\mu =  \max_{\substack{  v_{h} \in V_h|_{\Omega_h \setminus \mathcal{D}} \\  v_{h} |_{ \partial \mathcal{D} }= w_h }}
\frac{\int \limits_{\mathcal{D}} \left| \nabla w_h \right|^2 dx  + \sum \limits_{s=1}^m \frac{1}{d_s^2} \left[ \int \limits_{\mathcal{D}^s} w_h dx \right]^2}{
\int \limits_{\mathcal{D}} \left| \nabla w_h \right|^2 dx  + \int \limits_{\Omega_h \setminus \mathcal{D}} \left| \nabla v_{h} \right|^2 dx}
\geqslant 
\frac{\norm{w_h}_d^2}{\norm{w_h}_d^2  + \int \limits_{\Omega_h \setminus \mathcal{D}} \left| \nabla v_{h} \right|^2 dx},
\label{e:kuz_big_big_est}
\end{equation}
for any $v_{h} \in V_h |_{\Omega_h \setminus \mathcal{D}}$, such that $v_{h} = w_h$ on $ \partial \mathcal{D}$, where
\begin{equation}
 \norm{w_h}_d^2 = \int \limits_{\mathcal{D}} \left| \nabla w_h \right|^2 dx  + \sum \limits_{s=1}^m \frac{1}{d_s^2} \left[ \int \limits_{\mathcal{D}^s} w_h dx \right]^2. 
\end{equation}
\begin{figure}[h!tb]
	\centering
	\includegraphics[width=0.6\textwidth]{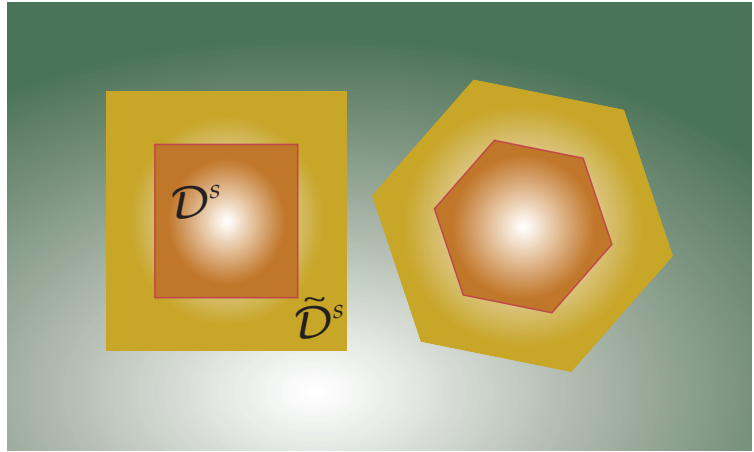}
	\caption{An example 	of $\mathcal{D}^s$ and $ \tilde{\mathcal{D}}^s$.}
	\label{f:kuz_fig_ex1}
\end{figure}
Following  \cite{kuz09,kuz18}, we embed subdomains $\mathcal{D}^s$ into subdomains $\tilde{\mathcal{D}}^s$ with the conforming boundary $\tilde{\Gamma}=\partial\tilde{\mathcal{D}}^s$ (see Figure \ref{f:kuz_fig_ex1}) so that
\begin{equation}\label{E:d_st_estimate}
 \min \limits_{x \in \ol{\mathcal{D}}^s, ~y \in \tilde{\mathcal{D}}^s  \cup  \Gamma} |x - y | \geqslant c d_s, 
\end{equation}
with a given positive constant $c$ independent of $d_s$, $s\in \{1,\ldots,m\}$. We assume that $ \tilde{\mathcal{D}}^s\cap \tilde{\mathcal{D}}^t = \varnothing $ for any $s \neq t$, $s,t \in \{1,\ldots,m\}$. We define  $\tilde{\mathcal{D}} = \bigcup \limits_{s=1}^m \tilde{\mathcal{D}}^s$, and assume that $v_{h}$ in (\ref{e:kuz_big_big_est}) vanishes in $\Omega_h \setminus \tilde{\mathcal{D}}$.
With that, we obtain the following estimate 
\begin{equation*}
\mu \geqslant  \min \limits_{s\in \{1,\ldots,m\}} \frac{\norm{w_{h,s}}_{d,s}^2}{\norm{w_{h,s}}_{d,s}^2+ \int \limits_{\tilde{\mathcal{D}}^s \setminus \mathcal{D}^s} | \nabla v_{h} |^2 dx},
\end{equation*}	
for any $v_{h} \in V_h|_{\tilde{\mathcal{D}}^s \setminus \mathcal{D}^s}$, such that  $v_{h}|_{ \Gamma_{s}} =  w_{h,s}$, \hspace{1pt} $v_{h} |_{ \tilde{\Gamma}_{s}}=  0$, where $w_{h,s}:=w_h|_{\mathcal{D}^{s}}$, and  
\begin{equation*}
 \norm{w_{h,s}}^2_{d,s}= \int \limits_{\mathcal{D}^{s}} | \nabla w_{h,s} |^2 dx + \frac{1}{d_s^2} \left[\int \limits_{\mathcal{D}^{s}} w_{h,s} dx \right]^2. 
\end{equation*}
If we assume that for any $w_{h,s} \in V_s$, its finite element extension $ \tilde{w}_{h} \in V_h|_{ \tilde{\mathcal{D}}^s  \setminus \mathcal{D}^{s}}$ with $\tilde{w}_{h}|_{\Gamma_{s} }= w_{h,s} $, $\tilde{w}_{h}|_{ \tilde{\Gamma}_s} = 0$, exists such that
\begin{equation}
\int \limits_{\tilde{\mathcal{D}}^s  \setminus \mathcal{D}^{s}} | \nabla \tilde{w}_{h}	|^2 dx \leqslant C^2 \norm{w_{h,s}}^2_{d,s}, 
\label{E:kuz_label_new_new_eq}
\end{equation} 
with a positive constant $C$ independent of $\Omega_h$ and values of $d_s$, $s\in \{1,\ldots,m\}$, then we arrive at the estimate
\begin{equation}\label{E:eig_est_mu_c}
\mu \geqslant \frac{1}{1+C^2}.
\end{equation}
The existence of  norm preserving finite element extensions on quasi-uniform regular shaped triangular meshes was proved in \cite{toswid05}.
To utilize the latter result to (\ref{E:kuz_label_new_new_eq})  we have to assume that the mesh $\Omega_h$ is quasi-uniform and regular shaped in subdomains $\tilde{\mathcal{D}}^s  \setminus \overline{\mathcal{D}}^{s} $ and to apply the transformation $x' = \frac{1}{d_s} x$ for each of the subdomains $\tilde{\mathcal{D}}^s $ as it was proposed in \cite{kuz09}, $s\in \{1,\ldots,m\}$.

Thus, under the assumptions made, the estimate
\begin{equation*}
a_0 \geqslant \frac{1}{1+C^2}
\end{equation*} 
holds, where $C$ is a positive constant independent on $\Omega_h$ and values of $d_s$, $s\in \{1,\ldots,m\}$.

\begin{remark}
There is an alternative proof of the estimate for $\mu$ from below as in \eqref{E:eig_est_mu_c}, see \cite{kuz18}, that does not use the algebraic technique \eqref{e:kuz_eigen_mass}-\eqref{e:kuz_equival_eigenval_pr} proposed in this paper.
\end{remark}

\subsection{Eigenvalue estimates for the matrix $\boldsymbol{\mathcal{H}} \boldsymbol{\mathcal{A}}_{\varepsilon}$}
\label{S:eig_est_ha}

In this Section, we assume that the assumptions made in the end of the Section \ref{S:eig_est_hs} are still valid, that is, the mesh $\Omega_h$ in $\tilde{\mathcal{D}}^s$ is regularly shaped and quasi-uniform, $s \in \{ 1, \ldots, m\}$, and distances  between ${\mathcal{D}}^s$ and ${\mathcal{D}}^t$ satisfy  (\ref{E:d_st_estimate}) with a constant $c$ independent of $\Omega_h$ as well as shape and location of inclusions. In other words, we assume that
\begin{equation}\label{E:eig_est_s0_bd}
\boldsymbol{\mathrm{S}}_0 \leqslant \boldsymbol{\mathcal{B}}_{\mathcal{D}} \leqslant \left(1 + C^2 \right) \boldsymbol{\mathrm{S}}_0,
\end{equation}
where $C^2$ is a positive constant independent of $\Omega_h$, and shape and locations of ${\mathcal{D}}^s$, $s \in \{ 1, \ldots, m\}$.
Consider the eigenvalue problem
\begin{equation}\label{E:eig_est_a_eps}
\boldsymbol{\mathcal{A}}_{\varepsilon} \begin{bmatrix}
\overline{v} \\
\overline{w}
\end{bmatrix}
=
\mu \boldsymbol{\mathcal{H}}^{-1}_0  \begin{bmatrix}
\overline{v} \\
\overline{w}
\end{bmatrix},
\end{equation}  
and two additional eigenvalue problems
\begin{equation}\label{E:eig_est_a_hat_eq}
\boldsymbol{\hat{\mathcal{A}}} \begin{bmatrix}
\overline{v} \\
\overline{w}
\end{bmatrix}
=
\hat{\mu} \boldsymbol{\mathcal{H}}^{-1}_0  \begin{bmatrix}
\overline{v} \\
\overline{w}
\end{bmatrix},
\end{equation}
\begin{equation}\label{E:eig_est_a_check_eq}
\boldsymbol{\check{\mathcal{A}}} \begin{bmatrix}
\overline{v} \\
\overline{w}
\end{bmatrix}
=
\check{\mu} \boldsymbol{\mathcal{H}}^{-1}_0  \begin{bmatrix}
\overline{v} \\
\overline{w}
\end{bmatrix},
\end{equation}  
with the matrices 
\begin{equation*}
\boldsymbol{\hat{\mathcal{A}}} = 
\begin{bmatrix} \bold{A} & \bold{B}^T \\ 
\bold{B} &  - \bold{Q}  
\end{bmatrix} ,
\quad
\mbox{and} 
\quad
\boldsymbol{\check{\mathcal{A}}} = 
\begin{bmatrix} \bold{A} & \bold{B}^T \\ 
\bold{B} & -r_{\max} \boldsymbol{\mathrm{S}}_0 - \bold{Q}  
\end{bmatrix} ,
\end{equation*}  
respectively, where
\begin{equation*}
r_{\max} = \left(1+C^2\right) \varepsilon_{\max} ,
\end{equation*}
and
\begin{equation*}
\boldsymbol{\mathcal{H}}^{-1}_0 = \begin{pmatrix}
\bold{A} & 0 \\0 & \bold{B} \bold{A}^{-1}  \bold{B}^T + \bold{Q}  
\end{pmatrix}.
\end{equation*}
It is obvious that
\begin{equation*}
\boldsymbol{\check{\mathcal{A}}} \leqslant\boldsymbol{\mathcal{A}}_\varepsilon \leqslant \boldsymbol{ \hat{\mathcal{A}}} ,
\end{equation*}
and three eigenproblems (\ref{E:eig_est_a_eps}), (\ref{E:eig_est_a_hat_eq}), and (\ref{E:eig_est_a_check_eq}) have equal numbers of negative and positive eigenvalues. It is also obvious that all three eigenproblems have the same multiplicity of the eigenvalue $\mu = \check{\mu}=\hat{\mu} = 1$ and the underlying eigenvectors $\begin{bmatrix}
\overline{v} \\
\overline{w}
\end{bmatrix} $ satisfy the conditions 
\begin{equation}\label{E:hat_mu_est}
\overline{v} \in \ker \bold{B}, \quad  \overline{w} \in \ker \bold{B}^T = \ker \boldsymbol{\mathcal{B}}_{\mathcal{D}}.
\end{equation}
The latter condition in  (\ref{E:hat_mu_est}) implies that for the eigenvalues  $\mu , \check{\mu},\hat{\mu}$ not equal to one in \eqref{E:eig_est_a_eps}-\eqref{E:eig_est_a_check_eq}, we can impose additional conditions on the  vector $\overline{w} \in \mathbb{R}^n$: 
\begin{equation}\label{E:eig_prob_add_cond}
\left(\bold{M}\overline{w} , \overline{\xi} \right) = 0 , \quad \forall \overline{\xi} \in \ker \boldsymbol{\mathcal{B}}_{\mathcal{D}}. 
\end{equation}
Assume that  $\hat{\mu} \neq 1$ in (\ref{E:eig_est_a_hat_eq}) and  $\check{\mu} \neq 1$ in (\ref{E:eig_est_a_check_eq}), then eliminating the vector $\overline{v} \in \mathbb{R}^N$ in (\ref{E:eig_est_a_hat_eq}) and (\ref{E:eig_est_a_check_eq}) (see also \cite{kuz95}) yields the equations
\begin{equation*}
-\frac{1}{1-\hat{\mu}} = \hat{\mu}, 
\quad \mbox{and} \quad
-\frac{1}{1- \check{\mu}} - r_{\max} =  \check{\mu} ,
\end{equation*}
respectively. It follows that each of eigenproblems (\ref{E:eig_est_a_hat_eq}) and (\ref{E:eig_est_a_check_eq}) under the condition \eqref{E:eig_prob_add_cond} has only two different eigenvalues
\begin{equation} \label{E:mu12} 
\hat{\mu}_{1,2} = \frac{1 \mp \sqrt{ 5}}{2},
\quad \mbox{and} \quad
\check{\mu}_{1,2} = \frac{1-r_{\max} \mp \sqrt{ \left(1-r_{\max}\right)^2 + 4 \left(1+r_{\max}\right)}}{2},
\end{equation}
respectively.

\begin{remark}
 It is obvious that   $\check{\mu}_1$ tends to $\hat{\mu}_1= \frac{1}{2} \left(1 - \sqrt{5}\right)$  and $\check{\mu}_2$ tends to $\hat{\mu}_2 =\frac{1}{2} \left(1 + \sqrt{5}\right)$ as  $\varepsilon_{\max}$ tends to zero. 
\end{remark}
Straightforward analysis of (\ref{E:mu12})  shows that
\begin{equation*}
\check{\mu}_1 < \hat{\mu}_1 < 0 < \check{\mu}_2 < \hat{\mu}_2.
\end{equation*}
Using inequalities \eqref{E:eig_est_s0_bd} and results of  \cite{bell}, we conclude that all eigenvalues of (\ref{E:eig_est_a_hat_eq}), which are not equal one, belong to the union of two disjoint segments
\begin{equation*}
\left[ \check{\mu}_1, \hat{\mu}_1 \right] \cup \left[\check{\mu}_2, \hat{\mu}_2 \right],
\end{equation*}
with the endpoints independent of $\Omega_h$, shape, and  location of the inclusions (see the assumption in the beginning of this Section).
Simple analysis shows that  $\check{\mu}_2 >1$ for any $\varepsilon_{\max} \geqslant 0$. To this end, we conclude that all the eigenvalues of (\ref{E:eig_est_a_eps}) belong to the set
\begin{equation*}
\left[\check{\mu}_1, \hat{\mu}_1\right] \cup \left[1,\hat{\mu}_2\right].
\end{equation*}
Now we consider the eigenvalue problem
\begin{equation}\label{E:eig_est_a_eps_main}
\boldsymbol{\mathcal{A}}_{\varepsilon} \begin{bmatrix}
\overline{v} \\
\overline{w}
\end{bmatrix}
=
\mu \boldsymbol{\mathcal{H}}^{-1}  \begin{bmatrix}
\overline{v} \\
\overline{w}
\end{bmatrix} ,
\end{equation}  
where 
\begin{equation}\label{E:h_a_defin}
\boldsymbol{\mathcal{H}}^{-1}=
\begin{bmatrix}
\mathcal{H}_{\mathrm{A}}^{-1} &0 \\ 0 & \boldsymbol{\mathcal{B}}_{\mathcal{D}} +  \bold {Q}
\end{bmatrix} ,
\end{equation}
and $\mathcal{H}_{\mathrm{A}}$ is a symmetric positive definite matrix satisfying the condition
\begin{equation} \label{E:beta-bounds}
\beta_1 \bold{A} \leqslant \mathcal{H}_{\mathrm{A}}^{-1} \leqslant \beta_2 \bold{A}, 
\end{equation}
with positive constants $\beta_1$ and $\beta_2$. We assume that $\beta_1$ and $\beta_2$ are independent of $\Omega_h$. For instance, $\mathcal{H}_{\mathrm{A}}^{-1}$ could be BPX or AMG preconditioner \cite{bpx90,bpx97,bs08,kuz90}.
Using (\ref{E:eig_est_s0_bd}) and (\ref{E:h_a_defin}), we obtain
\begin{equation*}
\alpha_{\min} \boldsymbol{\mathcal{H}}^{-1}_0 \leqslant \boldsymbol{\mathcal{H}}^{-1} \leqslant \alpha_{\max} \boldsymbol{\mathcal{H}}^{-1}_0,
\end{equation*}
where 
\begin{equation*}
\alpha_{\min} = \min \Big\{ \beta_1; \frac{1}{1+C^2} \Big\}, 
\quad \mbox{and} \quad
\alpha_{\max} = \max \{\beta_2;1\} .
\end{equation*}
Then, straightforward analysis shows that the eigenvalues of the matrix $\boldsymbol{\mathcal{H}} \boldsymbol{\mathcal{A}}_{\varepsilon}$ belong to the set
\begin{equation*}
\left[C_1, C_2\right] \cup \left[ C_3, C_4\right],
\end{equation*}  
where
\begin{equation} \label{E:constC}
C_1 = \frac{\check{\mu}_1}{\alpha_{\min}}\leqslant C_2 = \frac{\hat{\mu}_{\textcolor{red}{1}}}{\alpha_{\max}}< 0,
\quad \mbox{and} \quad
C_4 = \frac{\hat{\mu}_2}{\alpha_{\min}} > C_3 = \frac{1}{\alpha_{\max}}>0 .
\end{equation}
Thus, we have proved the following result.

\begin{theorem} \label{T:main}
Let the mesh $\Omega_h$ be regularly shaped and quasi-uniform, and distances  between ${\mathcal{D}}^s$ and ${\mathcal{D}}^t$ satisfy  (\ref{E:d_st_estimate}) with a constant $c$ independent of $\Omega_h$ as well as the shape and location of inclusions. 
Then the eigenvalues of the matrix  $\left(\boldsymbol{\mathcal{H}}  \boldsymbol{\mathcal{A}}_{\varepsilon}\right)^2$ belong to the segment   $\left[ a^2, b^2 \right]$, where
\begin{equation*}
a = \min \big\{|C_2|; C_3 \big\}, \quad  b = \max \big\{ |C_1|; C_4 \big\} ,
\end{equation*}
where $C_i$, $i\in \{1,\ldots,4\}$, are given by \eqref{E:constC}.
\end{theorem}
Note that under the assumption made,  the values of $a$ and $b$ are independent of $\Omega_h$ as well as shape and location of inclusions ${\mathcal{D}}^s$ in $\Omega$, $s\in \{1,\ldots,m\}$.

\begin{remark} 
The results of this section can be easily extended to the case of 3D diffusion problem as well as to the problems with nonzero reaction coefficient and to different types of boundary conditions (Neumann, Robin and mixed).
\end{remark} 

\section{Numerical results} \label{S:numerics}
	
\begin{figure}%
\centering
\subfigure[$m=16$, $h = \frac{1}{64}$, $ d= 8h$]{%
\label{f:kuz_fig_ex}%
\includegraphics[height=2.75in]{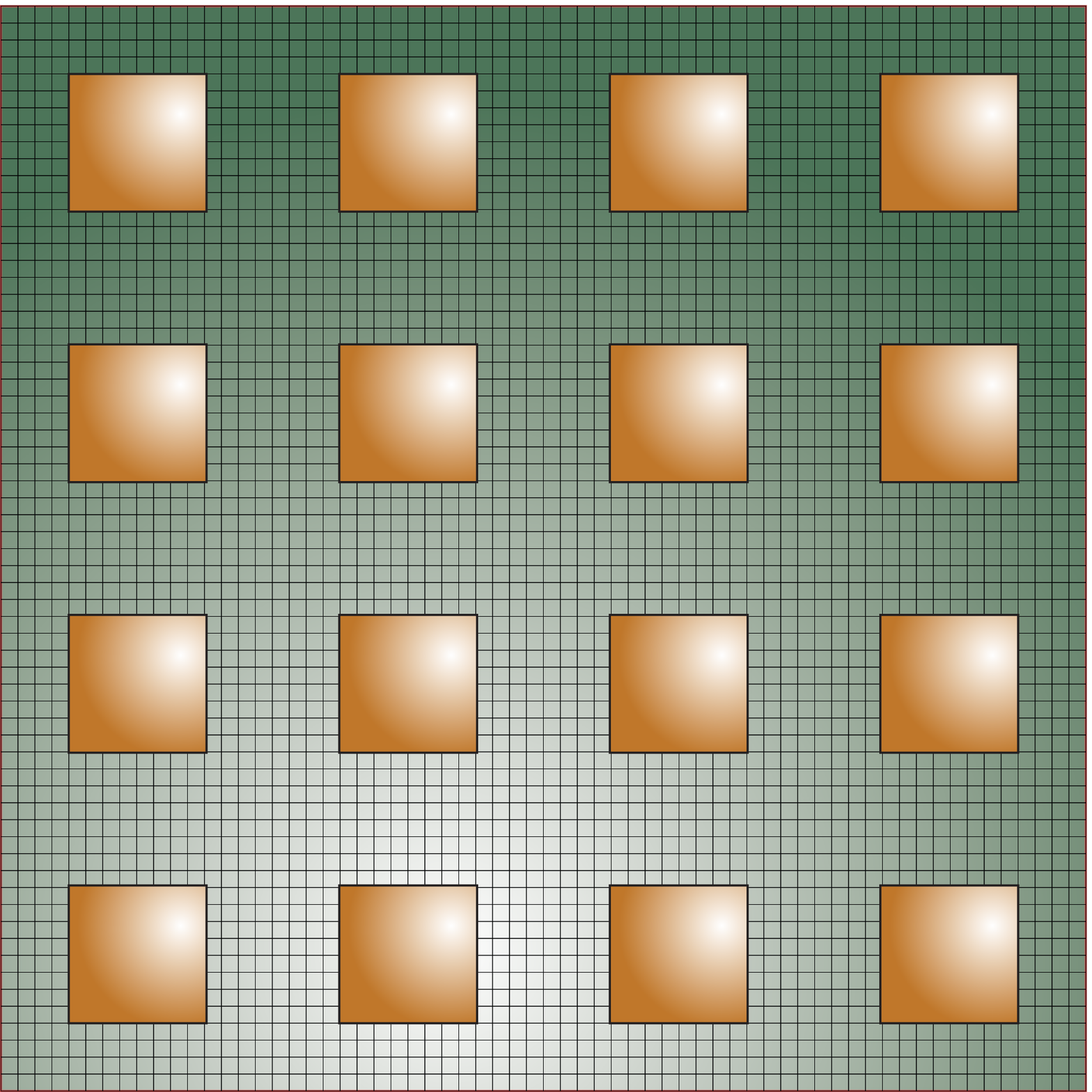}}%
\quad
\subfigure[$m=256$, $h = \frac{1}{64}$, $ d= 2h$]{%
\label{f:kuz_quad_reg}%
\includegraphics[height=2.75in]{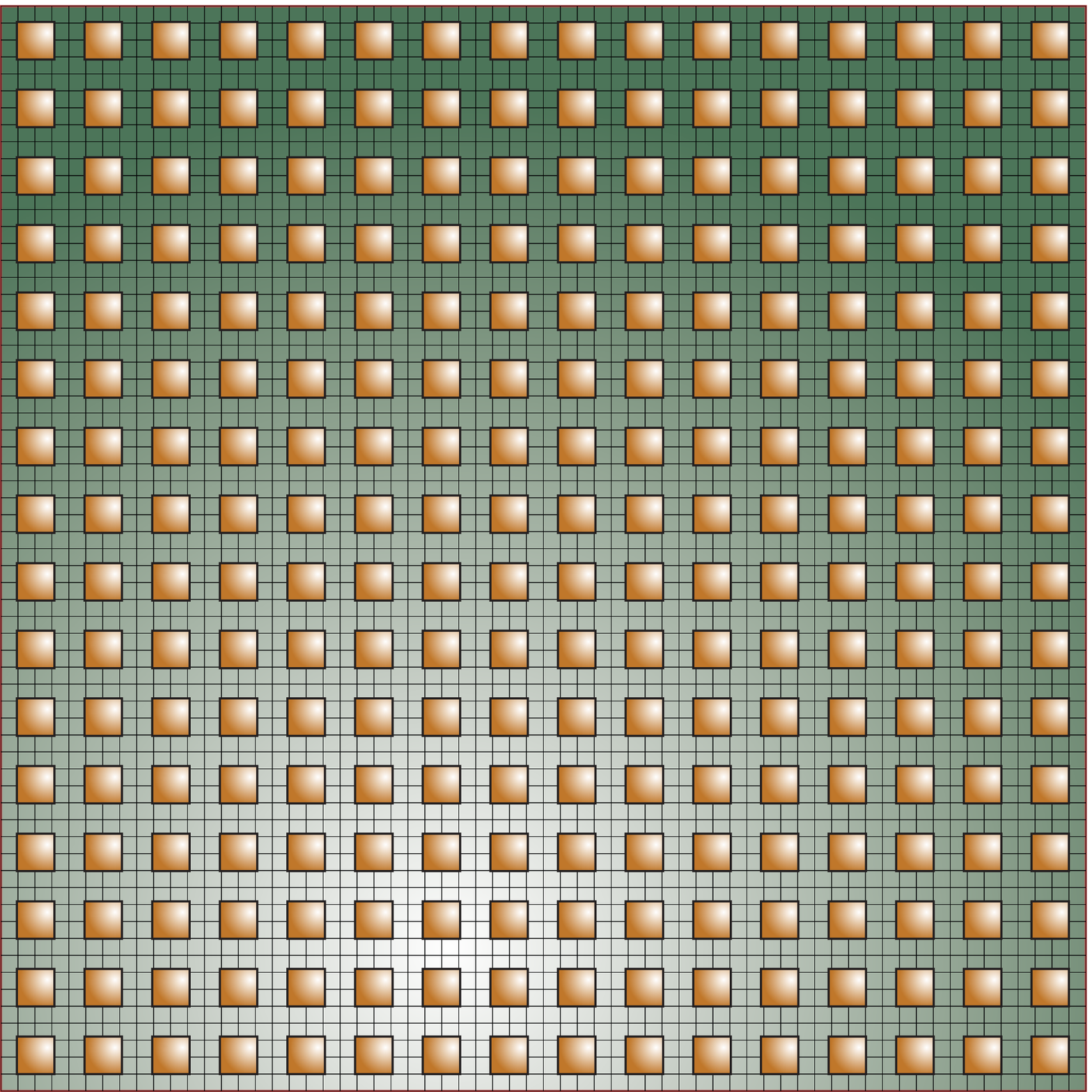}}%
\caption{Periodic distributions of particles }
\end{figure}

\begin{figure}[h!tb]
	\centering
	\includegraphics[width=0.45\textwidth]{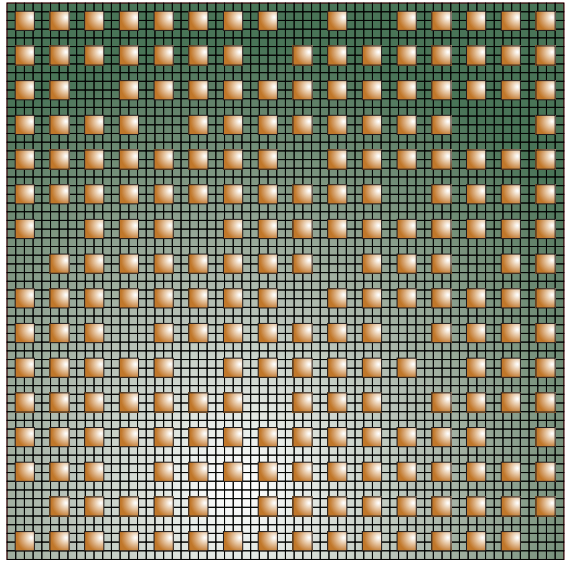}
	\caption{Random distribution of particles ($m=230$, $h = \frac{1}{64}$, $ d= 2h$)}
	\label{f:kuz_quad_rand}
\end{figure}

To evaluate and verify methods proposed in Sections \ref{S:form} and \ref{S:methods}, and the theoretical results justified in Section \ref{S:eig_est}, we consider the following simple model problem. Let $\Omega$ be a unit square, and $\Omega_h$ be a triangulated square mesh with mesh step size $h = \frac{1}{\sqrt{N} -1}$. We consider two types of particles' distribution in $\Omega$, $s \in \{1,\ldots,m \}$. The first one, called {\it ``periodic''}, is shown in Figures \ref{f:kuz_quad_reg} and \ref{f:kuz_fig_ex}.  The second one, called {\it ``random''}, is obtained by removing $\hat{m} < \bar{m}$ inclusions randomly chosen from the periodic array of $\bar{m}$ particles, so that $m=\bar{m}-\hat{m}$, see Figure \ref{f:kuz_quad_rand}. The values of $\varepsilon_s$ in  $\mathcal{D}^s$, $s \in \{1,\ldots,m \}$, are chosen either randomly from the segment $\left[\varepsilon_{\min},10^{-2} \right]$,  where $\varepsilon_{\min} <1$, or uniformly $\varepsilon_s = \varepsilon$, $s \in \{1,\ldots,m\}$.

\begin{table}[h!]
\centering
\begin{tabular}{|c |c| c |c |c |} 
 \hline 
\diagbox[width=3em]{ $\delta$ }{N} & 65,025 & 261,121 & 1,046,529 & 4,190,209\\
\hline
$10^{-2}$ & 4 & 4 & 4  &4 \\  [1ex]
$10^{-4}$ & 7 & 7 & 7  & 7  \\  [1ex]
$10^{-6}$ & 10 & 10 & 10  & 10  \\  [1ex]
$10^{-7}$ & 12 & 12 & 12  & 12 \\  [1ex]
$10^{-8}$ & 14 & 14 & 14  & 14  \\  [1ex]
\hline 
\end{tabular}
\caption{The number of PCG iterations}
\label{T:table:1}
\end{table}

In our numerical tests the inclusions are represented by $d \times d$ squares separated by the distance $d\equiv d_s$, $s \in \{1,\ldots,m \}$, between neighboring inclusions so that the minimal distance between the inclusions and the boundary $\partial \Omega$ equals $d/2$ as shown in Figure \ref{f:kuz_quad_reg}.

\begin{table}[h!]
\centering
\begin{tabular}{|c |c| c |c |c | c| c|} 
\hline 
\multirow{ 2}{*}{\diagbox[width=5em]{$\varepsilon_{\min}$}{$m$}} & \multicolumn{2}{c|}{65,536} & \multicolumn{2}{c|}{16,384} & \multicolumn{2}{c|}{4,096}\\
\cline{2-7}
& Period&  Rand&  Period&  Rand &  Period&  Rand \\  
\hline                                                
$10^{-2}$ & 11 & 11 & 11  &10 & 10 & 10  \\  [1ex]
$10^{-4}$ & 11 & 11 & 11  & 11& 10  & 10   \\  [1ex]
$10^{-6}$ & 11 & 11 & 11  & 11  & 10 & 10 \\  [1ex]
\hline 
\end{tabular}
\caption{The number of  PU iterations}
\label{T:table:2}
\end{table}

\begin{table}[h!]
\centering
\begin{tabular}{|c |c| c |c |c | c| c|} 
\hline 
\multirow{ 2}{*}{\diagbox[width=5em]{$\varepsilon_{\min}$}{$m$}} & \multicolumn{2}{c|}{65,536} & \multicolumn{2}{c|}{16,384} & \multicolumn{2}{c|}{4,096}\\
\cline{2-7}
&  Period&  Rand&  Period&  Rand &  Period&  Rand \\  
\hline                                                
$10^{-2}$ & 40 & 40 & 43  &43 & 46 & 44  \\  [1ex]
$10^{-4}$ & 40 & 40 & 44  & 44& 46  & 46  \\  [1ex]
$10^{-6}$ & 40 & 40 & 44  & 44  & 46 & 46 \\  [1ex]
\hline 
\end{tabular}
\caption{The number of PL iterations}
\label{T:table:3}
\end{table}

\begin{table}[h!]
\centering	
\begin{tabular}{|c |c| c |c |c | c| c|} 
\hline 
\multirow{ 2}{*}{\diagbox[width=5em]{$\varepsilon_{\min}$}{$m$}} & \multicolumn{2}{c|}{65,536} & \multicolumn{2}{c|}{16,384} & \multicolumn{2}{c|}{4,096}\\
\cline{2-7}
&  Period&  Rand&  Period&  Rand &  Period&  Rand \\  
\hline                                                
$10^{-2}$ & 90 &90 & 88  &88 &89 & 89  \\  [1ex]
$10^{-4}$ & 93 & 93 & 92  & 92& 92  & 92  \\  [1ex]
$10^{-6}$ & 93 & 93 & 92  & 92  & 92 & 92 \\  [1ex]
\hline 
\end{tabular}
\caption{The number of PCG iterations}
\label{T:table:4}
\end{table}

The matrix $\mathcal{H}_{\mathrm{A}}$ is the W-cycle Algebraic Multigrid preconditioner, proposed and investigated in \cite{kuz90,kuz_pde}. It was shown in \cite{kuz_pde}, that the eigenvalues of the matrix $\mathcal{H}_{\mathrm{A}} \bold{A}$ lie in the segment $\left[\frac{1}{2}\left(3-\sqrt{3}\right), \frac{3}{2}\left(1+\sqrt{3}\right)\right]$, that is, in \eqref{E:beta-bounds} we have
\begin{equation*}
\beta_1 = \frac{1}{2}\left(3-\sqrt{3}\right), \qquad \beta_2 = \frac{3}{2}\left(1+\sqrt{3}\right).
\end{equation*}
Therefore, the number of arithmetical operations (flops) for calculation of the matrix-vector product $\mathcal{H}_{\mathrm{A}} \ol{\xi}$ with $\ol{\xi} \in \mathbb{R}^N$ is bounded above by $5\times N$, hence, arithmetical costs of multiplication of a vector by $\mathcal{H}_{\mathrm{A}}$ and $\bold{A}$ are almost equal.

The main goal of our numerical experiments is to evaluate the minimal number of iterations sufficient for the minimization of initial errors in $\delta^{-1}$ times, $\delta<1$. To this end, in our numerical tests, we consider the homogeneous systems with a randomly chosen initial guess.

For the PU method (\ref{E:PCG_Uzawa-0})-(\ref{E:PCG_Uzawa-2}) the stopping criteria was 
\begin{equation} \label{E:stop-crit1}
\norm{\ol{p}^k}_{\bold{S}_\varepsilon} \leqslant \delta \norm{\ol{p}^0}_{\bold{S}_\varepsilon},
\end{equation}
and for the PL (\ref{E:Prec_lanc-0})-(\ref{E:Prec_lanc-3}) and PCG method (\ref{E:prec_pcg-0})-(\ref{E:prec_pcg-2}) the stopping criteria was 
\begin{equation}  \label{E:stop-crit2}
\norm{\ol{z}^k}_{\bold{K}_\varepsilon} \leqslant \delta \norm{\ol{z}^0}_{\bold{K}_\varepsilon} .
\end{equation}

In Table \ref{T:table:1}, we display the number of PCG iterations with the preconditioner $\mathcal{H}_\mathrm{A}$ mentioned at the beginning of this section for the homogeneous system
\begin{equation*}
\bold{A} \ol{x} = \ol{0} ,
\end{equation*}
and randomly chosen initial guesses $\ol{x}^0$. The stopping criteria was 
\begin{equation*}
\norm{\ol{x}^k}_{\bold{A}}\leqslant \delta \norm{\ol{x}^0}_{\bold{A}} .
\end{equation*}
We observe that $12$ iterations are sufficient to minimize the $\bold{A}$-norm of the error in $10^7$ times.

In Table \ref{T:table:2}, we display the number of iterations of the PU method with $\delta = 10^{-6}$, which is independent of a random choice of $\varepsilon \in \left[\varepsilon_{\min}, 10^{-2} \right]$ in the algebraic system, and the distribution of the inclusions. 
To perform the product $\mathcal{H}_{\mathrm{A}} \ol{\xi}$, $\ol{\xi} \in \mathbb{R}^N$, we used $12$ iterations of the PCG method for systems with the matrix $\bold{A}$.
\begin{table}[h!]
\centering
\begin{tabular}{|c |c| c |c |} 
\hline 
\multirow{ 2}{*}{\diagbox[width=6.5em]{$\varepsilon_{\min}$}{method}} &  & & \\
 & \textbf{PL} & PCG & PU\\
\cline{2-4}
\hline                                                
$10^{-2}$ & \textbf{44} &176 & 120    \\  [1ex]
$10^{-4}$ & \textbf{46} & 184 & 132   \\  [1ex]
$10^{-6}$ & \textbf{46} & 184 & 132   \\  [1ex]
\hline 
\end{tabular}
\caption{Arithmetical cost}
\label{T:table:5}
\end{table}

In Tables \ref{T:table:3} and  \ref{T:table:4}, we display the number of iterations for the PL and PCG methods described in Sections \ref{S:PL} and \ref{S:PCG}, respectively. The tests are done for various numbers of particles $m$, and the two types of particles' distribution: periodic and random ones. As it is clearly seen, the number of iterations does not depend on $\varepsilon_{\min} $, nor on distribution of the particles, or  their number, or the mesh size $h$ in $\Omega_h$.

Using results of the tests presented in Tables \ref{T:table:2},  \ref{T:table:3} and  \ref{T:table:4}, we compare all three respective methods (PU, PL, and PCG) in terms of their arithmetical costs  in Table \ref{T:table:5}. 
Note that due to Remark \ref{R:rk1}, the major computational effort is associated with multiplications by the matrices $\mathcal{H}_\mathrm{A}$ and $\bold{A}$, hence, this table presents the number of multiplications by $\mathcal{H}_\mathrm{A}$ and $\bold{A}$ needed to solve the underlying systems with accuracy $\delta$ due to criteria \eqref{E:stop-crit1} and  \eqref{E:stop-crit2}. 
Based on these results, we may conclude that for the above test problems, the PL method is almost three times faster than the PU method, and almost four times faster than the PCG method. Obviously, 
the results and conclusions may be different for other test problems and different choice of a preconditioner  $\mathcal{H}_{\mathrm{A}}$  for the matrix $\bold{A}$.

\section{Conclusions} \label{S:concl}
This paper proposes three preconditioned iterative methods for solving a linear system of the saddle point type arising in discretization of the diffusion problem   \eqref{E:pde-form} that involves large variation of its coefficient \eqref{E:sigma}. The latter feature is typically called {\it high contrast}.  
The main theoretical outcome presented in Theorem \ref{T:main} yields  that with the proposed preconditioner $\boldsymbol{\mathcal{H}}$, the
condition numbers of the preconditioned matrix $\boldsymbol{\mathcal{H}}\boldsymbol{\mathcal{A}}_\varepsilon$ is of $O(1)$. This implies robustness of the proposed preconditioners.
The assumption about regularly shaped and quasi-uniform mesh $\Omega_h$ is needed to apply the norm-preserving extension theorem of \cite{widlund87} that yields independence of convergence rates of the mesh size $h$. In order to claim independence of convergence rates of 
the diameter of $\mathcal{D}^s$, $s\in \{1,\ldots,m\}$ and their locations, we need assumption \eqref{E:d_st_estimate}.
Our numerical experiments based on simple test scenarios presented in Section \ref{S:numerics} confirm theoretical findings of this paper, and demonstrate convergence rates of the proposed iterative schemes to be independent of the contrast, discretization size, and also on the number of  inclusions and their sizes. The very important feature of the discussed procedures is that they are computationally inexpensive with the arithmetical cost being proportional to the size of  the linear system. This makes the proposed methodology attractive for the type of applications that use high contrast particles.

\end{document}